\documentclass[12pt,twoside]{article}
\usepackage[frame,cmtip,arrow,matrix,line,graph,curve]{xy}
\usepackage{amssymb}

\date{}

\begin{document}

\centerline{\bf }

\centerline{\bf }

\centerline{}

\centerline{}

\centerline{\Large{\bf On Pseudo-Umbilical Spacelike Submanifolds in }}

\centerline{}

\centerline{\Large{\bf  Indefinite Space Form $M_{p}^{n+p}(c)$}}

\centerline{}

\centerline{\bf {Majid Ali Choudhary}}

\centerline{}

\centerline{Department of Mathematics}

\centerline{Zakir Husain Delhi College(Evening), New Delhi -110025 (India)   }

\centerline{E-mail : majid\_alichoudhary@yahoo.co.in}

\centerline{}

\newtheorem{Theorem}{\quad Theorem}[section]

\newtheorem{Definition}[Theorem]{\quad Definition}

\newtheorem{Corollary}[Theorem]{\quad Corollary}

\newtheorem{Lemma}[Theorem]{\quad Lemma}

\newtheorem{Example}[Theorem]{\quad Example}

\centerline{}

\begin{abstract} In the present note, first we derive an intrinsic inequality for Pseudo-umbilical spacelike submanifold in an indefinite space form. We use this inequality to show that such submanifold is totally geodesic. In the rest part of this paper, using a result of Aiyama \cite{1}, we prove that Pseudo-umbilical spacelike subamnifold is totally umbilical.

\end{abstract}



\section{Introduction}
\setcounter{equation}{0}
\renewcommand{\theequation}{1.\arabic{equation}}
\hspace {1cm}

Let $M^n$ be an $n$-dimensional Riemannian manifold immersed in an $(n+p)$-dimensional connected semi-Riemannian manifold $M_p^{n+p}(c)$ of constant curvature $c$ whose index is $p$. We call $M_p^{n+p}(c)$ a space form of index $p$ and simply a space form when $p=0$. If $c > 0$, we call it a de Sitter space of index $p$. As the semi-Riemannian metric of $M_p^{n+p}(c)$ induces the Riemannian metric of $M^n$, $M^n$ is called a spacelike submanifold. Let $h$ be the second fundamental form of the immersion and $N$ be the mean curvature vector. Denote by $<.,.>$  the scalar product of $M^n$. If there exists a function $\lambda\geq0$  on $M^n$  such that
 \begin{equation}
    <h(X,Y),N>=\lambda <X,Y>        \label{1.1}
 \end{equation}
for any tangent vectors $X,Y$ on $M^n$, then $M^n$ is called a Pseudo-umbilical spacelike submanifold of $M_p^{n+p}(c)$. If the mean curvature vector $N$ vanishes identically, then $M^n$  is called a maximal spacelike submanifold of $M_p^{n+p}(c)$. Every maximal spacelike submanifold of $M_p^{n+p}(c)$ is itself a Pseudo-umbilical spacelike submanifold of $M_p^{n+p}(c)$.

Spacelike hypersurfaces and submanifolds have attracted the attention of many mathematicians in the recent years e.g. Dong [7], Wu ([15],[16]), Liu [9]. In the year 2002, Pang [13] studied spacelike hypersurfaces in de sitter space and derived an intrinsic inequality to obtain a sufficient and necessary condition for such hypersurfaces to be totally geodesic. Han [8], investigated spacelike submanifolds in indefinite space form $M_p^{n+p}(c)$ and obtained an intrinsic inequality to prove some rigidity theorem. However, Pseudo-umbilical submanifolds have also been paid attention by many mathematicians e.g. [2], [3], [14].   X. F. Cao [2], gave an intrinsic inequality for pseudo umbilical spacelike submanifolds in the indefinite space form. In 1995, Sun [14] first proved that mean curvature $H$ of the pseudo-umbilical submanifolds in indefinite space form $M_p^{n+p}(c)$ is constant.
On the other hand, Y. Zheng [17] gave an intrinsic condition for a compact space like hypersurface in a de sitter space to be totally umbilical. While, Ximin ([10],[11],[12]) extended Cheng-Yau [6] technique to investigate spacelike hypersurfaces and spacelike submanifolds with constant scalar curvature and proved some intrinsic conditions for such hypersurface or submanifold to be totally umbilical.

Inspired by all the above investigations, in the first half of this paper, we give an intrinsic inequality for pseudo-umbilical spacelike submanifold of indefinite space form $M_p^{n+p}(c)$. Using this inequality, we get a necessary and sufficient condition for such a submanifold to be totally geodesic. In the rest part of this note, Using Chen-Yau [6] technique and taking into account the results obtained by Aiyama [1] and Sun [14], we prove that pseudo-umbilical spacelike submanifold of indefinite space form with constant scalar curvature and flat normal bundle is totally umbilical.

\section{Preliminaries}
\setcounter{equation}{0}
\renewcommand{\theequation}{2.\arabic{equation}}

We choose a local field of semi-Riemannian orthonormal frame $e_1,...,e_{n+p}$ in $M_p^{n+p}(c)$ such that at each point of $M^n$, $e_1,...,e_n$ span the tangent space of $M^n$ and form an orthonormal frame there. We use the following convention on the range of indices:\\
$1\leq A,B,C,...\leq n+p$;\hspace {1cm}  $1\leq i,j,k,...\leq n$; \hspace {1cm} $n+1\leq \alpha,\beta,\gamma\leq n+p.$\\
Let $\omega_1,...,\omega_{n+p}$ be its dual frame field so that the semi-Riemannian metric of $M_p^{n+p}(c)$ is given by
\begin{eqnarray*}
d\bar{s}^2=\sum_{i}\omega^2_{i}-\sum_{\alpha}(\omega_{\alpha})^2=\sum_{A}\epsilon_{A}\omega^2_{A} \hspace {.5cm},
\end{eqnarray*}
where $\epsilon_i=1$ and $\epsilon_\alpha=-1$.

Then the structure equations of $M_p^{n+p}(c)$ can be written as:
\begin{eqnarray*}
d\omega_A=\sum_{B}\epsilon_B\omega_{AB}\wedge\omega_B \hspace {.5cm}, \hspace {.5cm} \omega_{AB}+\omega_{BA}\\
d\omega_{AB}=\sum_{C}\epsilon_C\omega_{AC}\wedge\omega_{CB}-\frac{1}{2}\sum_{C,D}K_{ABCD}\omega_C\wedge\omega_D\\
K_{ABCD}=c \epsilon_A \epsilon_B(\delta_{AC}\delta_{BD}-\delta_{AD}\delta_{BC})
\end{eqnarray*}
Restricting these forms to $M^n$, we obtain $\omega_\alpha=0$, \hspace {.2cm} $n+1\leq\alpha\leq n+p$ \hspace {.2cm} and the Riemannian metric of $M^n$ is written as
$ ds^2=\sum_{i}\omega^2_{i} $.
  From Cartan's   lemma we can write \begin{eqnarray*} \omega_{\alpha i}=\sum_{j}h^\alpha_{ij}\omega_j, \hspace {.3cm} h^\alpha_{ij}=h^\alpha_{ji} \end{eqnarray*}
From these formulas, the structure equations of $M^n$ are given by:
\begin{eqnarray}
d\omega_i=\sum_{j}\omega_{ij}\wedge\omega_j \hspace {.5cm}, \hspace {.5cm} \omega_{ij}+\omega_{ji}\nonumber\\
d\omega_{ij}=\sum_{k}\omega_{ik}\wedge\omega_{kj}-\frac{1}{2}\sum_{k,l}K_{ijkl}\omega_k\wedge\omega_l\nonumber\\
R_{ijkl}=c(\delta_{ik}\delta_{jl}-\delta_{il}\delta_{jk})-\sum_{\alpha}(h^\alpha_{ik}h^\alpha_{jl}-h^\alpha_{il}h^\alpha_{jk})
\end{eqnarray}
$R_{ijkl}$ , being the components of the curvature tensor of $M^n$.\\
The second fundamental form of $M^n$ is given by
\begin{eqnarray*}
h=\sum_{i,j,\alpha}h^\alpha_{ij}\omega_i\otimes\omega_j\otimes e_\alpha
\end{eqnarray*}
The mean curvature vector $N$ of $M^n$  is defined by
\begin{eqnarray*}
N=\frac{1}{n}\sum trH_\alpha e_\alpha.
\end{eqnarray*}
Here \textit{tr} is the trace of the matrix $H_\alpha=(h^\alpha_{ij})$ and it is well known that $N$ is independent of the choice of unit normal vectors $e_1,...,e_{n+p}$ to $M^n$. The length of the mean curvature vector is called the mean curvature of $M^n$  and is denoted by $H$.
Now, let $e_{n+1}$ be parallel to $N$. Then we have
\begin{eqnarray}
trH_{n+1}=nH, \hspace {.3cm} trH_\alpha=0, \hspace {.3cm} \alpha\neq0.
\end{eqnarray}
Define the first and second covariant derivatives of ${h^\alpha_{ij}}$, say ${h^\alpha_{ijk}}$ and ${h^\alpha_{ijkl}}$ by
\begin{eqnarray*}
\sum_{k}h^\alpha_{ijk}\omega_{k}=dh^\alpha_{ij}+\sum_{k}h^\alpha_{kj}\omega_{ki}+\sum_{k}h^\alpha_{ik}\omega_{kj}+\sum_{\beta}h^\beta_{ij}\omega_{\beta\alpha}
\end{eqnarray*}
\begin{eqnarray*}
\sum_{l}h^\alpha_{ijkl}\omega_{l}&=&dh^\alpha_{ijk}+\sum_{m}h^\alpha_{mjk}\omega_{mi}+\sum_{m}h^\alpha_{imk}\omega_{mj}\\
&&+\sum_{m}h^\alpha_{ijm}\omega_{mk}+\sum_{\beta}h^\beta_{ijk}\omega_{\beta\alpha}
\end{eqnarray*}
Then we have
\begin{eqnarray}
h^\alpha_{ijk}&=&h^\alpha_{ikj}\\
h^\alpha_{ijkl}-h^\alpha_{ijlk}&=&\sum_{m}h^\alpha_{im}R_{mjkl}+\sum_{m}h^\alpha_{jm}R_{mikl}+\sum_{\beta}h^\beta_{ij}R_{\alpha\beta kl}
\end{eqnarray}
where $R_{\alpha\beta kl}$ are the components of the normal curvature tensor of $M^n$ , that is $R_{\alpha\beta kl}=\sum_{i}(h^\alpha_{ik}h^\beta_{il}-h^\alpha_{il}h^\beta_{ik})$. If $R_{\alpha\beta kl}=0$ at a point $x$ of $M^n$ , we say that normal connection of $M^n$  is flat at $x$, and it is well known that $R_{\alpha\beta kl}=0$ at $x$ if and only if $h^\alpha_{ij}$ are simultaneously diagonalizable at $x$. [4]

The Laplacian $\Delta h^\alpha_{ij}$ of the fundamental form $h^\alpha_{ij}$  is defined to be $\sum_{k}h^\alpha_{ijkk}$ and hence, if $M^n$ has flat normal bundle, then from (2.3) and (2.4), we have
\begin{eqnarray}
\Delta h^\alpha_{ij}&=&\sum_{k}(h^\alpha_{ijkk}-h^\alpha_{ikjk})+\sum_{k}(h^\alpha_{ikjk}-h^\alpha_{ikkj})+\sum_{k}(h^\alpha_{ikkj}-h^\alpha_{kkij})+\sum_{k}h^\alpha_{kkij}\nonumber\\ &=&\sum_{m,k} h^\alpha_{im}R_{mkjk}+\sum_{m,k} h^\alpha_{mk}R_{mijk}+\sum_{k}h^\alpha_{kkij} .
\end{eqnarray}
But, in view of (1.1) and (2.2), we have
\begin{eqnarray*}\langle h(e_i,e_j),H_{e_{n+1}}\rangle=H^2\delta_{ij} \end{eqnarray*}
so, we arrive at the following
 \begin{eqnarray*}\sum h^\alpha_{ij}h^\alpha_{kkij}=nH\Delta H.	\end{eqnarray*}	
As H is constant due to Sun [14], we get
\begin{eqnarray}\sum h^\alpha_{ij}h^\alpha_{kkij}=0.	\end{eqnarray}
Since, normal bundle of $M^n$ is flat, we can diagonalize the second fundamental form simultaneously, so that $h^\alpha_{ij}=\lambda^\alpha_{i}\delta_{ij}, \alpha=n+1,...,n+p$ and then using (2.5), we have
\begin{eqnarray}
h^\alpha_{ij}\Delta h^\alpha_{ij}=\frac{1}{2}\sum R_{ijij}(\lambda^\alpha_i-\lambda^\alpha_j)^2 .
\end{eqnarray}

\section{Pseudo-umbilical spacelike submanifold}
\setcounter{equation}{0}
\renewcommand{\theequation}{3.\arabic{equation}}

In order to prove our result, we state the following lemma.

\begin{Lemma} Let $a_1,...,a_n$ be real numbers, then
\begin{eqnarray*}\sum (a_i)^2\geq \frac{1}{n}(\sum a_i)^2 \end{eqnarray*}

and the equality holds if and only if $a_1=...=a_n$.
\end{Lemma}

We prove the follwing.

\begin{Theorem} Let $M^n$ be n-dimensional compact Pseudo-umbilical spacelike submanifold in $M^{n+p}_p(c)  (c > 0)$, $S$ and $\rho$ be Ricci curvature tensor and scalar curvature of $M^n$, respectively, then
 \begin{eqnarray*}|S|^2\geq2c\rho(n-1)-c^2n(n-1)^2.\end{eqnarray*}
\end{Theorem}
{\bf Proof.} From the Gauss equation, we derive [2, eq. (2.7)]
\begin{eqnarray*}S_{ij}=(n-1)c\delta_{ij}+\sum_{k,\alpha}h^\alpha_{ik}h^\alpha_{jk}\end{eqnarray*}
So, we 	have	
\begin{eqnarray*}|S|^2&=&\sum S^2_{ij}\\
&=&\sum_{ij}\{(n-1)c\delta_{ij}+\sum_{k,\alpha}h^\alpha_{ik}h^\alpha_{jk}\}^2\\
&=&n(n-1)^2c^2+\sum_{ij}(\sum_{k,\alpha}h^\alpha_{ik}h^\alpha_{jk})^2+2c(n-1)(\sum_{i,k,\alpha}h^\alpha_{ik}h^\alpha_{ik})
\end{eqnarray*}
and the scalar curvature is given by
\begin{eqnarray*}\rho&=&\sum_i S_{ii}\\
&=&cn(n-1)+\sum_{ij\alpha}(h^\alpha_{ij})^2
\end{eqnarray*}
above equations can be rewritten in the following way
\begin{eqnarray*}|S|^2&=&n(n-1)^2c^2+\sum_{ij}(\sum_{k,\alpha}h^\alpha_{ik}h^\alpha_{jk})^2+2c(n-1)(\rho-cn(n-1))\\
&=&2c(n-1)\rho-n(n-1)^2c^2+\sum_{ij}(\sum_{k,\alpha}h^\alpha_{ik}h^\alpha_{jk})^2
\end{eqnarray*}                  		
Since $M^n$ has flat normal bundle, we can diagonalize the second fundamental form simultaneously, so that $h^\alpha_{ij}=\lambda^\alpha_{i}\delta_{ij}$ and then we have using lemma 3.1
\begin{eqnarray*}|S|^2&=&2c(n-1)\rho-n(n-1)^2c^2+\sum_{ij}(\sum_{k,\alpha}\lambda^\alpha_{i}\lambda^\alpha_{j}\delta_{ik}\delta_{jk})^2\\
&=&2c(n-1)\rho-n(n-1)^2c^2+\sum_i(\sum_{\alpha}\lambda^{\alpha^2}_i)^2\\
&\geq&2c(n-1)\rho-n(n-1)^2c^2+\frac{1}{n}(\sum_i(\sum_{\alpha}\lambda^{\alpha^2}_i))^2
\end{eqnarray*}or, we can write the above equation as follows
\begin{eqnarray}|S|^2&=&2c(n-1)\rho-nc^2(n-1)^2+\frac{1}{n}(\sum_{i,\alpha}\lambda^{\alpha^2}_i)^2 \nonumber\\
|S|^2&\geq&2c(n-1)\rho-nc^2(n-1)^2 \end{eqnarray}                  		
whereby proving the result.
\hfill $\Box$

Next, we prove
\begin{Theorem} Let $M^n$ be Pseudo-umbilical spacelike submanifold in $M^{n+P}_p (c)$ and $S$ and $\rho$ be Ricci curvature tensor and scalar curvature of $M^n$  , respectively, then $|S|^2=2c(n-1)\rho-nc^2(n-1)^2$  if and only if $M^n$ is totally geodesic.
\end{Theorem}
{\bf Proof.} If $M^n$ is totally geodesic, that is $\lambda^{\alpha}_i=0$,\hspace {.3cm} $i \in \{1,...,n\}$ then, we have
\begin{eqnarray*}|S|^2=nc^2(n-1)^2 	\hspace {.3cm} and \hspace {.3cm}  \rho=nc(n-1)\end{eqnarray*}
such that $|S|^2=2c(n-1)\rho-nc^2(n-1)^2$.
Conversely, if equality holds in (3.1), then all the inequalities of (3.1) become equality. From lemma 3.1, we have
\begin{eqnarray}(\sum_{i,\alpha}\lambda^{\alpha^2}_i)^2=0 \end{eqnarray} and
\begin{eqnarray} \lambda^{\alpha^2}_1=\lambda^{\alpha^2}_2=...=\lambda^{\alpha^2}_n \end{eqnarray}                  		
for $i,j \in \{1,...,n\}$ and $\alpha \in \{n+1,...,n+p\}$.

In the light of (3.2) and (3.3), we conclude that $\sum\lambda^{\alpha^2}_i=0$, which shows that	 $\lambda^{\alpha}_i=0$, whereby proving that $M^n$ is totally geodesic.
\hfill $\Box$

\section{Pseudo-umbilical spacelike submanifold with constant mean curvature}
\setcounter{equation}{0}
\renewcommand{\theequation}{4.\arabic{equation}}

First we prove the following lemma which shall be used later to prove the main result.

\begin{Lemma} Let $M^n$ be n-dimensional compact Pseudo-umbilical spacelike submanifold in $M^{n+p}_p(c)$ with mean curvature $H$. If normalized scalar curvature $\rho$ is constant and $\rho<c$, then
\begin{eqnarray} \sum_{i,j,k}h^2_{ijk}\geq n^2|\nabla H|^2 \end{eqnarray}
\end{Lemma}
{\bf Proof.} Using equation (2.1), we can easily see that
\begin{eqnarray*}n^2H^2-\|h\|^2=n(n-1)(c-\rho). \end{eqnarray*}
Taking the covariant derivative of above equation and using the fact that $\rho$ is constant, we obtain
\begin{eqnarray*}n^2HH_k=\sum_{i,j,\alpha}h^\alpha_{ij}h^\alpha_{ijk} \hspace {.3cm} k=1,...,n \end{eqnarray*}	
and hence using Cauchy-Schwartz inequality, we have
\begin{eqnarray}\sum_k n^4H^2H^2_k=(\sum_{i,j,\alpha}h^\alpha_{ij}h^\alpha_{ijk})^2 \leq |h|^2 \sum_{i,j,\alpha}(h^\alpha_{ijk})^2 \hspace {.3cm} \end{eqnarray}	 where equality holds if and only if there exists a real function $c_k$ such that
\begin{eqnarray*} h^\alpha_{ijk}=c_k h^\alpha_{ij} \end{eqnarray*}for all $i,j$ and $\alpha$. Taking sum on both sides of (4.2) with respect to $k$, we get
\begin{eqnarray}n^4H^2|\nabla H|^2= n^4H^2\sum_k H^2_k\leq |h|^2 \sum_{i,j,k,\alpha}(h^\alpha_{ijk})^2\leq n^2H^2\sum_{i,j,k,\alpha}(h^\alpha_{ijk})^2.\end{eqnarray}Therefore, (4.1) holds on $M^n$.
\hfill $\Box$

Now, we prove the main result.

\begin{Theorem}Let $M^n$ be $n$-dimensional compact Pseudo-umbilical spacelike submanifold with mean curvature $H$ immersed in $M^{n+p}_p(c)$. Suppose that $M^n$   has flat normal bundle and scalar curvature $\rho$ is constant and $\rho<c$, then $M^n$  is totally umbilical and isometric to a sphere.
\end{Theorem}
\textbf{Proof. }Since, the Laplacian of $|h|^2$  is given by
\begin{eqnarray*} \frac{1}{2}|h|^2=\sum_{\alpha,i,j,k}(h^\alpha_{ijk})^2+\sum_{\alpha,i,j}h^\alpha_{ij}\triangle	 h^\alpha_{ij}\end{eqnarray*}
So, in the light of (2.7), above equation reduces to
\begin{eqnarray} \frac{1}{2}|h|^2=\frac{1}{2}n^2\triangle H^2=\|\nabla h\|^2+\frac{1}{2}\sum_{i,j,\alpha}R_{ijij}(\lambda^\alpha_i-\lambda^\alpha_j)^2\end{eqnarray}
Now, define an operator $\Box$ acting on $f$ by
\begin{eqnarray*} \Box f=\sum_{ij}(nH\delta_{ij}-h^{n+1}_{ij})f_{ij}\end{eqnarray*}
Since, $\sum_{ij}(nH\delta_{ij}-h^{n+1}_{ij})$  is trace free it follows from [6] that the operator $\Box$ is self adjoint relative to  $L^2$-inner product of    $M^n$, that is
\begin{eqnarray*} \int f\Box g=\int g\Box f \end{eqnarray*}
Thus, we have

\begin{eqnarray}
\Box H&=&\sum_{ij}(nH\delta_{ij}-h^{n+1}_{ij})H_{ij}\nonumber\\ &=&nH\sum_i H_{ii}-\sum_i \lambda^{n+1}_i H_{ii}\\ &=&\frac{1}{2}n(\Delta H^2-2\|\nabla H\|^2)-\sum_i \lambda^{n+1}_i H_{ii}\nonumber\end{eqnarray}
Now, taking account of equations (4.4) and (4.5), we have
\begin{eqnarray*}
\Box(nH)=|\nabla h|^2-n^2|\nabla H|^2+\frac{1}{2}\sum R_{ijij}(\lambda^\alpha_i-\lambda^\alpha_j)^2\end{eqnarray*}
which on using (4.3) reduces to
\begin{eqnarray*}
\Box (nH)&\geq&\|h\|^{-2} n^3 H^2 |\nabla H|^2-n^2|\nabla H|^2+\frac{1}{2}\sum R_{ijij}(\lambda^\alpha_i-\lambda^\alpha_j)^2\\ &\geq&(\|h\|^{-2} nH^2-1)n^2 |\nabla H|^2+\frac{1}{2}\sum R_{ijij}(\lambda^\alpha_i-\lambda^\alpha_j)^2\end{eqnarray*}

Using the fact that $\Box$ is self adjoint, we conclude that

\begin{eqnarray*}
0\geq\int_{M^n}\{(\|h\|^{-2} nH^2-1)n^2 |\nabla H|^2+\frac{1}{2}\sum R_{ijij}(\lambda^\alpha_i-\lambda^\alpha_j)^2\}\end{eqnarray*}

But, [14] we have that $H$ is constant. Therefore, our result follows immediately from a result of Aiyama ([1], Theorem 3) and this completes the proof of our theorem. \hfill$\Box$

\centerline{}
\centerline{}

{\bf Acknowledgements.} The author is thankful to Department of Science and Technology, Government of India, for its financial assistance provided through Inspire Fellowship No. DST/INSPIRE Fellowship/2009/[xxv] to carry out this research work.


\begin{thebibliography}{99}

\bibitem{1} {Aiyama R., Compact space like submanifolds in a Pseudo-Riemannian sphere $S^{n+p}_{p}(c)$
            , \em Tokyo J. Math.} {\bf 18} (1995), 81--90.

\bibitem{2}{Cao X. F., Pseudo-umbilical spacelike submanifolds in an indefinite space form,
            \em Balkan Journal of Geometry and its Applications,} {\bf 6, No. 2} (2001), 117-121.
.
\bibitem{3}{Cao X. F., Pseudo-umbilical submanifolds of constant curvature Riemannian
            manifolds, \em Glasgow Math. J.,}{\bf  43} (2001), 129-133.

\bibitem{4}{Chen B. Y., Geometry of submanifolds, \em Marcel Dekker, New York,} (1973).

\bibitem{5}{Cheng Q. M., Complete spacelike submanifolds in a de Sitter space with parallel mean
            curvature vector,\em  Math. Z.,} {\bf 206}(1991), 333-339.

\bibitem{6}{Cheng S. Y., Yau  S. T., Hypersurfaces with constant scalar curvature,
            \em Math. Ann.,} {\bf 225} (1977), 195-204.

\bibitem{7}{Dong Y. X.,  Bernstein theorems for spacelike graphs with parallel mean curvature
            and controlled growth, \em J. Geom. Phys.,} {\bf  58}(2008), 324-333.

\bibitem{8}{Han Y., Spacelike submanifolds in indefinite space form $M^{n+p}_p (c)$,
          \em Archivum Mathematicum (Brno), Tomus,} {\bf  46} (2010), 79-86.

\bibitem{9}{Liu X. M.,  spacelike submanifolds with constant scalar curvature, \em C.R.Acad. Sci.
            Paris Ser. I Math.,}{\bf  332}(2001), 729-734.

\bibitem{13}{Pang H. D., Xu  S. L. ,Dai sh.,  Spacelike hypersurfaces in de sitter space,
            \em J. Geom. Phys. } {\bf 42} (2002), 78-84.

\bibitem{14}{Sun H., On spacelike submanifolds of a pseudo-Riemannian space form,
            \em Note Mat.,}{\bf 15, No. 2} (1995), 215-224.

\bibitem{15}{Wu B. Y., On the volume Gauss map image of spacelike submanifolds in de sitter
	space form, \em  J. Geom. Phys.,}  {\bf 53}(2005), 336-344.

\bibitem{16}{Wu B. Y., On the mean curvature of spacelike submanifolds in semi Riemannian
	manifolds, \em J. Geom. Phys.,} {\bf 56}(2006), 1728-1735.

\bibitem{10}{Ximin Liu, Spacelike subamnifolds with constant scalar curvature in the de sitter spaces,
            \em  J. Korean Math. Soc.,} {\bf  38}  (2001) No. 1, 135-146.

\bibitem{11}{Ximin Liu, Complete spacelike hypersurfaces with constant scalar curvature,
	        \em manuscripta math.,} {\bf 105} (2001), 367-377.

\bibitem{12}{Ximin Liu, Spacelike subamnifolds in de  sitter spaces, \em J. Phys. A: Math. Gen.,} {\bf  34}
	(2001) 5463-5468.

\bibitem{17}{Zheng Y., Space like hypersurfaces with constant scalar curvature in the de sitter
	          spaces, \em Differential Geometry and its Applications,} {\bf 6}(1996) 51-54.


\end{thebibliography}
\end{document}